# Nordhaus-Gaddum-type theorem for rainbow connection number of graphs

Lily Chen, Xueliang Li, Huishu Lian

Center for Combinatorics, LPMC

Nankai University, Tianjin 300071, P. R. China

Email: lily60612@126.com, lxl@nankai.edu.cn, lhs6803@126.com

**Abstract**

An edge-colored graph $G$ is rainbow connected if any two vertices are connected by a path whose edges have distinct colors. The rainbow connection number of $G$, denoted $rc(G)$, is the minimum number of colors that are used to make $G$ rainbow connected. In this paper we give a Nordhaus-Gaddum-type result for the rainbow connection number. We prove that if $G$ and $\overline{G}$ are both connected, then $4 \leq rc(G) + rc(\overline{G}) \leq n+2$. Examples are given to show that the upper bound is sharp for all $n \geq 4$, and the lower bound is sharp for all $n \geq 8$. For the rest small $n = 4, 5, 6, 7$, we also give the sharp bounds.

## 1 Introduction

All graphs considered in this paper are simple, finite and undirected. Undefined terminology and notations can be found in [1]. Let $G$ be a nontrivial connected graph with an edge coloring $c : E(G) \to \{1, 2, \cdots, k\}$, $k \in \mathbb{N}$, where adjacent edges may be colored the same. A path $P$ of $G$ is a *rainbow path* if no two edges of $P$ are colored the same. The graph $G$ is called *rainbow-connected* if for any two vertices $u$ and $v$ of $G$, there is a rainbow $u - v$ path. The minimum number of colors for which there is an edge coloring of $G$ such that $G$ is rainbow connected is called the *rainbow connection number*, denoted by $rc(G)$. Clearly, if a graph is rainbow connected, then it is also connected. Conversely, any connected graph has a trivial edge coloring that makes it rainbow connected, just by assigning each edge a distinct color. An easy observation is that if $G$ has $n$ vertices then $rc(G) \leq n-1$, since one may color the edges of a spanning tree with distinct colors, and color the remaining edges with one of the colors already used. It is easy to see that if $H$ is a connected spanning subgraph of $G$, then $rc(G) \leq rc(H)$. It is easy to see that $rc(G) = 1$ if and only if $G$ is a clique, and $rc(G) = n-1$ if and only if $G$ is a tree, as well as that a cycle with $k > 3$ vertices has a rainbow connection number $\lceil k/2 \rceil$. Also notice that $rc(G) \geq diam(G)$, where $diam(G)$ denotes the diameter of $G$.

A Nordhaus–Gaddum-type result is a (tight) lower or upper bound on the sum or product of the values of a parameter for a graph and its complement. The name "Nordhaus–Gaddum-type" is so given because it is Nordhaus and Gaddum [3] who first established the following type of inequalities for chromatic number of graphs in 1956. They proved that if $G$ and $\overline{G}$ are complementary graphs on $n$ vertices whose chromatic numbers are $\chi(G)$, $\chi(\overline{G})$, respectively, then

$$2\sqrt{n} \leq \chi(G) + \chi(\overline{G}) \leq n+1.$$

Since then, many analogous inequalities of other graph parameters are concerned, such as diameter [4], domination number [5], Wiener index and some other chemical indices [6], and so on. In this paper, we are concerned with analogous inequalities involving the rainbow connection number of graphs, we prove that

$$4 \leq rc(G) + rc(\overline{G}) \leq n+2.$$

The rest of this paper is organized as follows. First, we give the upper bound, and show that it is sharp for all $n \geq 4$. Then we give the lower bound, and show that it is also sharp for $n \geq 8$. Finally, for the rest small $n = 4, 5, 6, 7$, we give the sharp bound, respectively.

## 2  Upper bound on $rc(G) + rc(\overline{G})$

We know that if $G$ is a connected graph with $n$ vertices, then the number of the edges in $G$ must be at least $n-1$. So if both $G$ and $\overline{G}$ are connected then $n$ is not less than 4, since

$$2(n-1) \leq e(G) + e(\overline{G}) = e(K_n) = \frac{n(n-1)}{2}. \tag{$*$}$$

In the rest of the paper, we always assume that all graphs have at least 4 vertices, and both $G$ and $\overline{G}$ are connected.

**Lemma 1** $rc(G) + rc(\overline{G}) \leq n + 2$ for $n = 4, 5$, and the bound is sharp.

*Proof.* Note that $rc(G) \leq n-1$, equality holds if and only if $G$ is a tree. So

$$rc(G) + rc(\overline{G}) \leq 2(n-1),$$

equality holds if and only if both $G$ and $\overline{G}$ are trees. Then $(*)$ must holds with equality. That is, $n$ have to be 4, and

$$rc(G) + rc(\overline{G}) = 2(n-1) = 6 = 4 + 2.$$

Then

$$rc(G) + rc(\overline{G}) \leq 2n - 3$$

for $n \geq 5$.

For $n = 5$, let $G$ be a tree obtained from $S_4$ by attaching a pendent edge to one of the vertices of degree one. Then $rc(G) = 4$. We observe that $diam(\overline{G}) = 3$ and it can be colored by three colorings to make it rainbow connected. Thus $rc(\overline{G}) = 3$. Therefore, we have

$$rc(G) + rc(\overline{G}) = 7 = 2n - 3 = 5 + 2.$$

**Lemma 2** *Let $G$ be a nontrivial connected graph of order $n$, and $rc(G) = k$. Let $c : E(G) \to \{1, 2, \cdots, k\}$ be a rainbow $k$-coloring of $G$. Add a new vertex $P$ to $G$, $P$ is adjacent to $q$ vertices of $G$, the resulting graph is denoted by $G'$. Then if $q \geq n + 1 - k$, we have $rc(G') \leq k$.*

*Proof.* Let $X = \{x_1, x_2, \cdots, x_q\}$ be the vertices adjacent to $P$, $V \backslash X = \{y_1, y_2, \cdots, y_{n-q}\}$. If $q \geq n + 1 - k$, $n - q \leq k - 1$.

Since $G$ rainbow connected under the coloring $c$, for any $y_i$, $i \in \{1, 2, \cdots, n - q\}$, there is a rainbow $x_1 - y_i$ path, say $P_{x_1y_1}, P_{x_1y_2}, \cdots, P_{x_1y_{n-q}}$. For each $P_{x_1y_i}$, we find out the last vertex on the path that belongs to $X$, and the subpath between this vertex to $y_i$ of $P_{x_1y_i}$ is denoted by $P_i$. Then $P_i$ is a rainbow path whose vertices are in $Y$ except the first vertex.

Let $G_{x_i}$ be the union of the paths in $P_1, P_2, \cdots, P_{n-q}$ whose origin vertex is $x_i$, $1 \leq i \leq q$. If there is no path with origin vertex $x_i$, let $G_{x_i}$ be a trivial graph with the vertex $x_i$. Then $G_{x_i}$ is a subgraph of $G$, and $v(G_{x_i}) \leq n - q + 1 \leq k$. First, we consider the subgraph $G_{x_1}$, and let $V(G_{x_1}) = \{x_1, y_{i_1}, y_{i_2}, \cdots, y_{i_l}\}$.

**Case 1:** The number of colors appeared in $G_{x_1}$ is $k$. Then $e(G_{x_1}) \geq k$.

**Subcase 1.1:** $e(G_{x_1}) = k \geq v(G_{x_1})$.

In this case, $G_{x_1}$ contains a cycle, and no two edges of $G_{x_1}$ are colored the same. Thus, $G_{x_1}$ is rainbow connected. Let $e$ be an edge in the cycle. Then, by deleting $e$ and coloring the edge $Px_1$ with the color $c(e)$, we have that, for any $j \in \{1, 2, \cdots, l\}$, there is a rainbow $P - y_{i_j}$ path.

**Subcase 1.2** $e(G_{x_1}) > k \geq v(G_{x_1})$.

In this case, $G_{x_1}$ contains a cycle, and there are two edges $e_1, e_2$ with $c(e_1) = c(e_2)$.

If one of the edges, say $e_1$, is contained in a cycle of $G_{x_1}$. Then, by deleting it, we obtain a spanning subgraph $G'_{x_1}$ of $G_{x_1}$ with the same number of colors appearing in it, but $e(G'_{x_1}) = e(G_{x_1}) - 1$. If $e(G'_{x_1}) = k$, by a similar operation as in Subcase 1.1, we can obtain a coloring of $Px_1$ such that for any $j \in \{1, 2, \cdots, l\}$, there is a rainbow $P - y_{i_j}$ path. If $e(G'_{x_1}) > k$, we consider the graph $G'_{x_1}$ other than $G_{x_1}$.

If both $e_1, e_2$ are not in a cycle, they must be cut edges of $G_{x_1}$. Then, contract one of them, say $e_1$, and denote the resultant graph by $G'''_{x_1}$. The number of colors appeared in $G_{x_1}$ is still $k$, and $v(G'''_{x_1}) = v(G_{x_1}) - 1, e(G'''_{x_1}) = e(G_{x_1}) - 1$. If $e(G'''_{x_1}) = k$, by a similar operation as in Subcase 1.1, we can obtain a coloring of $Px_1$ such that for any $y_k$ in $G'''_{x_1}$, there is a rainbow $P - y_k$ path. It is easy to check that there still exists a rainbow $P - y_{i_j}$ path in $G_{x_1}$ for any $j \in \{1, 2, \cdots, l\}$. If $e(G'''_{x_1}) > k$, we consider the graph $G'''_{x_1}$ other than $G_{x_1}$.

**Case 2:** The number of colors appeared in $G_{x_1}$ is less than $k$. Then we color the edge $Px_1$ with a color not appeared in $G_{x_1}$.

No matter which cases happen, we can always color the edge $Px_1$ with one of the colors $\{1, 2, \cdots, k\}$, such that for any $j \in \{1, 2, \cdots, l\}$, there is a rainbow $P - y_{i_j}$ path.

For $G_{x_2}, G_{x_3}, \cdots, G_{x_q}$, we use the same way to color the edges $Px_2, Px_3, \cdots, Px_q$. Then we get a $k$-coloring of $G'$. Since for each $y_i$, there is an $x_j$, such that $y_i \in G_{x_j}$. Then the path $Px_jP_i$ is a rainbow path connecting $P$ and $y_i$. Thus in this coloring, $G'$ is rainbow connected. Therefore $rc(G') \leq k$.

**Theorem 1** $rc(G) + rc(\overline{G}) \leq n + 2$ for all $n \geq 4$, and this bound is best possible.

*Proof.* We use induction on $n$. From Lemma 1, the result is true for $n = 4, 5$. We assume that $rc(G) + rc(\overline{G}) \leq n + 2$ holds for complementary graphs on $n$ vertices. To the union of connected graphs $G$ and $\overline{G}$, a complete graph on the $n$ vertices, we adjoin a new vertex $P$. Let $q$ be the number of vertices of $G$ which are adjacent to $P$, then the number of vertices of $\overline{G}$ which are adjacent to $P$ is $n - q$. If $G'$ and $\overline{G}'$ are the resultant graphs (each of order $n + 1$), then

$$rc(G') \leq rc(G) + 1, rc(\overline{G}') \leq rc(\overline{G}) + 1.$$

These inequalities are evident from the fact that if given a rainbow $rc(G)$-coloring ($rc(\overline{G})$-coloring) of $G$ ($\overline{G}$), we assign a new color to the edges added from $P$ to $G$ ($\overline{G}$), the resulting coloring makes $G'$ ($\overline{G'}$) rainbow connected. Then $rc(G') + rc(\overline{G'}) \leq rc(G) + rc(\overline{G}) + 2 \leq n + 4$. And $rc(G') + rc(\overline{G'}) \leq n + 3$ except possibly when

$$rc(G') = rc(G) + 1, rc(\overline{G'}) = rc(\overline{G}) + 1.$$

In this case, by Lemma 2, $q \leq n - rc(G), n - q \leq n - rc(\overline{G})$, thus $rc(G) + rc(\overline{G}) \leq n$, from which $rc(G') + rc(\overline{G'}) \leq n + 2$. This completes the induction.

To see the bound can be attained, let $G$ be a tree obtained by joining the centers of two stars $S_p$ and $S_q$ by an edge $uv$, where $u$ and $v$ are the centers of $S_p$ and $S_q$, and $p + q = n$. Then $rc(G) = n - 1$. To compute the rainbow connection number of the complement graph of $G$, we assume that $X = V(S_p \backslash u) \cup \{v\}$, $Y = V(S_q \backslash v) \cup \{u\}$. Then $G$ is a bipartite graph with bipartition $(X, Y)$. Thus $\overline{G}[X]$, $\overline{G}[Y]$ is complete. We assign color 1 to $\overline{G}[X]$, 2 to $\overline{G}[Y]$ and 3 to the edges between $X$ and $Y$. The resulting coloring makes $\overline{G}$ rainbow connected, thus $rc(\overline{G}) \leq 3$. On the other hand $diam(G) = d(u, v) = 3$, it follows that $rc(\overline{G}) = 3$. Then we have $rc(G') + rc(\overline{G'}) = n + 2$.

# 3 Lower bound on $rc(G) + rc(\overline{G})$

As we have noted that $rc(G) = 1$ if and only if $G$ is a complete graph. In this case, $\overline{G}$ is not connected. Thus if both $G$ and $\overline{G}$ are connected, $rc(G) \geq 2$, $rc(\overline{G}) \geq 2$. That is, $rc(G) + rc(\overline{G}) \geq 4$.

**Proposition 1** *Let $G$ and $\overline{G}$ be complementary connected graphs with $rc(G) = rc(\overline{G}) = 2$. Then*
*(1) $diam(G) = diam(\overline{G}) = 2$.*
*(2) $2 \leq \delta(G) \leq \Delta(G) \leq n - 3$, $2 \leq \delta(\overline{G}) \leq \Delta(\overline{G}) \leq n - 3$.*
*(3) A vertex $u$ in $N_1(v)$ can not be adjacent to all vertices of $N_2(v)$, where $N_1(v)$, $N_2(v)$ is the first and second neighborhood of a vertex $v$, respectively.*

*Proof.* Since $2 \leq diam(G) \leq rc(G) = 2$, (1) clearly holds.

For (2), first, $\Delta(G) \neq n - 1$, otherwise $\overline{G}$ is disconnected.
Second, $\delta(G) \neq 1$. Indeed, if $\delta(G) = 1$, let $v$ be a vertex of degree one, and $u$ the vertex adjacent to $v$. Since $diam(G) = 2$, $u$ must be adjacent to all the other vertices, thus $d(u) = n - 1$, a contradiction. Similarly, $\delta(\overline{G}) \neq 1$. That is, $\delta(\overline{G}) \geq 2$. Therefore, $\Delta(G) \leq n - 1 - \delta(\overline{G}) \leq n - 3$, so does $\Delta(\overline{G})$.

For (3), if $u$ is adjacent to all vertices of $N_2(v)$, $u$ is not adjacent to them in $\overline{G}$, then the distance between $u$ and $N_2(v)$ is at least 2 in $\overline{G}$, and $v$ is adjacent to all vertices of $N_2(v)$, but not to the vertices in $N_1(v)$ of $\overline{G}$. So $d_{\overline{G}}(u, v) \geq 3$, which contradicts (1).

**Theorem 2** *For $4 \leq n \leq 7$, there are no graphs $G$ and $\overline{G}$ on $n$ vertices, such that $rc(G) = rc(\overline{G}) = 2$.*

*Proof.* We consider $n = 4, 5, 6, 7$, respectively.

**Case 1:** n=4.

Then, there is only one pair of complementary connected graphs, each is isomorphic to $P_4$, and its rainbow connection number is 3.

**Case 2:** n=5.

If $rc(G) = rc(\overline{G}) = 2$, by Proposition 1, $2 \leq \delta(G) \leq \Delta(G) \leq n - 3 = 2$. Then, $G \cong C_5$, $\overline{G} \cong C_5$. Since $rc(C_5) = 3$, there are no graphs $G$ and $\overline{G}$ on 5 vertices, such that $rc(G) = rc(\overline{G}) = 2$.

**Case 3:** n=6.

By Proposition 1, $2 \leq \delta(G) \leq \Delta(G) \leq n - 3 = 3$, the possible degree sequences are:

(a) $\begin{cases} d_G = (2,2,2,2,2,2) \\ d_{\overline{G}} = (3,3,3,3,3,3). \end{cases}$

(b) $\begin{cases} d_G = (3,3,2,2,2,2) \\ d_{\overline{G}} = (2,2,3,3,3,3). \end{cases}$

The graph $G$ with the degree sequence in (a) is a cycle of length 6, whose rainbow connection is 3. And the graph $G$ with the degree sequence in (b) satisfying Proposition 1 has to be the graph shown in Figure 1:

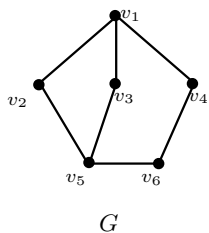

$G$

Figure 1: Graphs with degree sequence (3,3,2,2,2,2) satisfying Proposition 1.

Consider the pair of vertices $(v_2, v_4)$. The only 2-path is $v_2 v_1 v_4$, thus $c(v_1 v_2) \neq c(v_1 v_4)$. Similarly, $c(v_1 v_3) \neq c(v_1 v_4)$, then $c(v_1 v_2) = c(v_1 v_3)$. If we consider the pairs of vertices $(v_2, v_6)$, $(v_3, v_6)$, we have $c(v_2 v_5) = c(v_3 v_5)$. But then there is no rainbow $v_2 - v_3$ path, therefore $rc(G) \neq 2$.

**Case 4:** n=7.

By Proposition 1, $2 \leq \delta(G) \leq \Delta(G) \leq n - 3 = 4$, the possible degree sequences are: (in the following argument, we use two colors to color the edges of the graphs)

(1) $\begin{cases} d_G = (4, 4, 4, 4, 4, 4, 4) \\ d_{\overline{G}} = (2, 2, 2, 2, 2, 2, 2). \end{cases}$

In this case, $\overline{G}$ is a cycle of length 7, $rc(\overline{G}) = 4$.

(2) $\begin{cases} d_G = (4, 4, 4, 4, 4, 3, 3) \\ d_{\overline{G}} = (2, 2, 2, 2, 2, 3, 3). \end{cases}$

The graphs with the degree sequence (4,4,4,4,4,3,3) satisfying Proposition 1 are $G_1, G_2$ shown in Figure 2. The distance between $v_2$ and $v_5$ in $\overline{G_1}, \overline{G_2}$ is larger than 2, thus $rc(\overline{G_1}) \neq 2, rc(\overline{G_2}) \neq 2$.

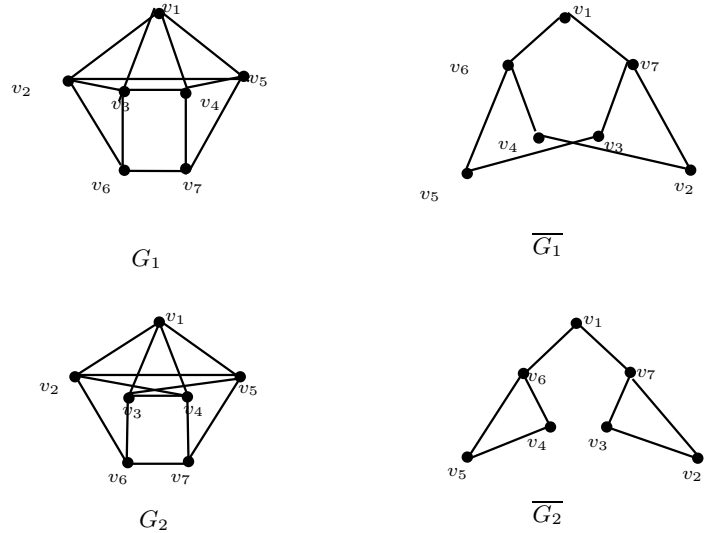

Figure 2: Graphs with degree sequence (4,4,4,4,4,3,3) satisfying Proposition 1.

(3) $\begin{cases} d_G = (4, 4, 4, 3, 3, 3, 3) \\ d_{\overline{G}} = (2, 2, 2, 3, 3, 3, 3). \end{cases}$

The graphs with degree sequence (4,4,4,3,3,3,3) satisfying Proposition 1 are subgraphs $G'_1, G'_2$ of $G_1, G_2$ shown in Figure 2 by deleting the edge $v_2 v_5$. We observe that the distance between $v_3$ and $v_4$ in $\overline{G'_1}, \overline{G'_2}$ is larger than 2, thus $rc(\overline{G'_1}) \neq 2, rc(\overline{G'_2}) \neq 2$.

(4) $\begin{cases} d_G = (4, 3, 3, 3, 3, 3, 3) \\ d_{\overline{G}} = (2, 3, 3, 3, 3, 3, 3). \end{cases}$

The graphs $G$ with degree sequence (4,3,3,3,3,3,3) satisfying Proposition 1 are $G_1, G_2$ and $G_3$ shown in Figure 3.

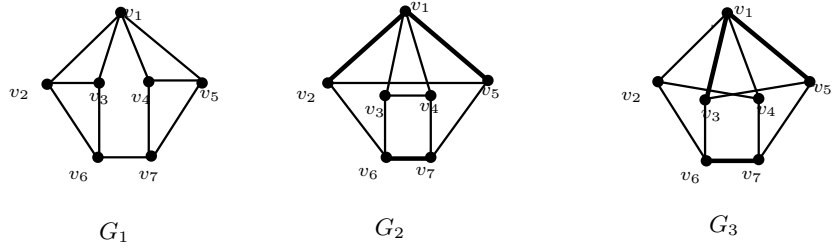

Figure 3: Graphs with degree sequence (4,3,3,3,3,3,3) satisfying Proposition 1.

Consider $G_1$. Since the only 2-path between $v_3$ and $v_4$ is $v_3v_1v_4$, $c(v_3v_1) \neq c(v_1v_4)$. Similarly, $c(v_3v_1) \neq c(v_1v_5)$. Then $c(v_1v_4) = c(v_1v_5)$, say color 2. By the same way, $c(v_3v_1) = c(v_1v_2) = 1$. Consider the pairs of vertices $(v_2, v_7), (v_3, v_7), (v_4, v_6), (v_5, v_6)$, we have $c(v_2v_6) = c(v_3v_6) = c(v_4v_7) = c(v_5v_7)$. If $c(v_2v_6) = 1$, there is no rainbow $v_1 - v_6$ path, and if $c(v_2v_6) = 2$, there is no rainbow $v_1 - v_7$ path. Therefore $rc(G_1) \neq 2$.

Consider $G_2$, whose rainbow connection number is 2, where the heavy lines is colored by color 2, the others are colored by color 1. So we consider its complement graph, the only 2-path between $v_6$ and $v_7$ is $v_6v_1v_7$, then $c(v_6v_1) \neq c(v_1v_7)$. Let $c(v_6v_1) = 1, c(v_1v_7) = 2$, thus $c(v_3v_7) = c(v_2v_7) = 1$, $c(v_4v_6) = c(v_5v_6) = 2$. If $c(v_2v_4) = 2$, there is no rainbow $v_2 - v_6$ path, and if $c(v_2v_4) = 1$, there is no rainbow $v_4 - v_7$ path. Therefore, we cannot use two colors to make $\overline{G_2}$ rainbow connected, that is $rc(\overline{G_2}) \neq 2$.

For $G_3$, whose rainbow connection number is also 2, by coloring the heavy lines with color 2, and the others with color 1. We consider its complement graph. By the same reason as above for $\overline{G_2}$, let $c(v_6v_1) = 1, c(v_1v_7) = 2, c(v_3v_7) = c(v_2v_7) = 1, c(v_4v_6) = c(v_5v_6) = 2$. Then, if $c(v_2v_5) = 2$, there is no rainbow $v_2 - v_6$ path, and if $c(v_2v_5) = 1$, there is no rainbow $v_5 - v_7$ path. Therefore, we cannot use two colors to make $\overline{G_2}$ rainbow connected, that is $rc(\overline{G_2}) \neq 2$.

(5) $\begin{cases} d_G = (4, 4, 4, 4, 3, 3, 2) \\ d_{\overline{G}} = (2, 2, 2, 2, 3, 3, 4). \end{cases}$

The graphs $G$ with degree sequence (4,4,4,4,3,3,2) satisfying Proposition 1 are $G_1, G_2, G_3, G_4$ shown in Figure 4.

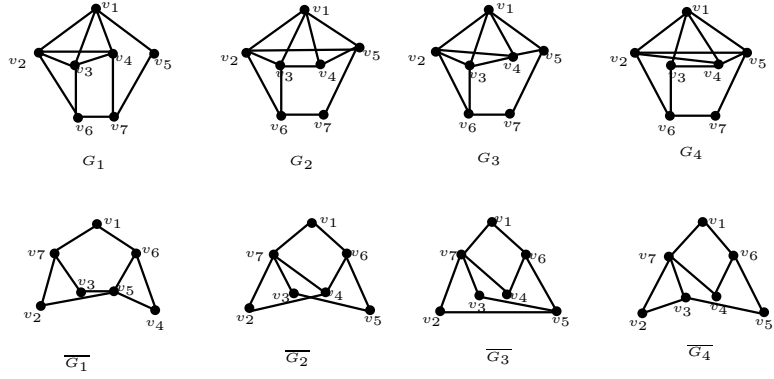

Figure 4: Graphs with degree sequence (4,4,4,4,3,3,2) satisfying Proposition 1.

Since $d_{\overline{G_1}}(v_4, v_7) = 3$, $d_{\overline{G_2}}(v_2, v_5) = 3$, $d_{\overline{G_4}}(v_2, v_6) = 3$, $rc(\overline{G_1}) \geq 3$, $rc(\overline{G_2}) \geq 3$, $rc(\overline{G_4}) \geq 3$.

For $\overline{G_3}$, consider the pair of vertices $(v_4, v_5)$. the only 2-path is $v_4 v_6 v_5$, so $c(v_4 v_6) \neq c(v_5 v_6)$, and let $c(v_5 v_6) = 2$. Similarly, $c(v_3 v_5) = c(v_2 v_5) = 1$, $c(v_7 v_2) = c(v_7 v_3) = 2$, $c(v_7 v_1) = c(v_7 v_4) = 1$,. If $c(v_1 v_6) = 2$, there is no rainbow $v_1 - v_5$ path, and if $c(v_1 v_6) = 1$, there is no rainbow $v_1 - v_4$ path. Thus, $rc(\overline{G_3}) \neq 2$.

(6) $\begin{cases} d_G = (4, 4, 3, 3, 3, 3, 2) \\ d_{\overline{G}} = (2, 2, 3, 3, 3, 3, 4). \end{cases}$

The graphs $\overline{G}$ with degree sequence (2,2,3,3,3,3,4) satisfying Proposition 1 have to be the following three graphs.

For $G_1$, consider the pair of vertices $(v_2, v_7)$. There is only one 2-path $v_2 v_6 v_7$ between them, then let $c(v_2 v_6) = 1, c(v_6 v_7) = 2$. Similarly, consider the pair of vertices $(v_3, v_7)$, we have $c(v_3 v_6) = 1$. Consider the pairs of vertices $(v_5, v_6)$, $(v_4, v_7)$, $(v_4, v_6)$, we get $c(v_5 v_7) = 1$, $c(v_4 v_5) = 2$, $c(v_3 v_4) = 2$. Consider the pairs of vertices $(v_2, v_5), (v_3, v_5)$, then $c(v_2 v_1) = c(v_3 v_1)$, thus there is no rainbow $v_2 - v_3$ path.

For $G_2$, consider the pairs of vertices $(v_4, v_6), (v_5, v_6)$, we have $c(v_4 v_7) = c(v_5 v_7)$, consider the pairs of vertices $(v_3, v_4), (v_3, v_5)$, we get $c(v_1 v_4) = c(v_1 v_5)$, thus there is no rainbow $v_4 - v_5$ path.

For $G_3$, its rainbow connection number is 2 by coloring the heavy lines with color 2, and assign color 1 to the other edges. So we consider its complement graph shown in the figure too. Since the only 2-path between $v_6$ and $v_7$ is $v_6 v_1 v_7$, let $c(v_1 v_7) = 1, c(v_1 v_6) = 2$. Thus $c(v_7 v_2) = c(v_7 v_3) = 2$, $c(v_6 v_4) = c(v_6 v_5) = 1$, $c(v_3 v_5) = 2$, $c(v_2 v_4) = 1$. If $c(v_2 v_5) = 1$, there is no rainbow $v_2 - v_6$ path. If $c(v_2 v_5) = 2$, there is no rainbow $v_5 - v_7$ path.

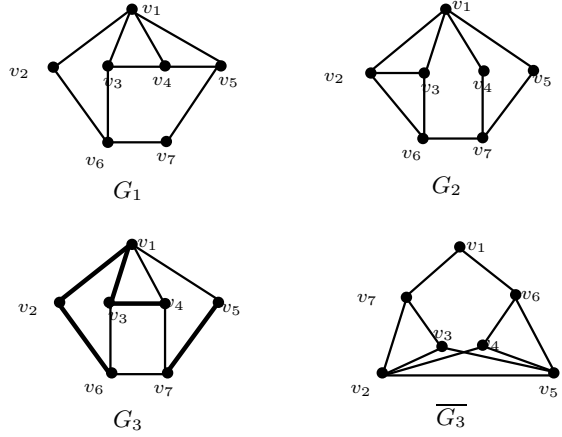

Figure 5: Graphs with degree sequence (2,2,3,3,3,3,4) satisfying Proposition 1.

(7) $\begin{cases} d_G = (4, 4, 4, 3, 3, 2, 2) \\ d_{\overline{G}} = (2, 2, 2, 3, 3, 4, 4). \end{cases}$

The graphs $\overline{G}$ with degree sequence (2,2,2,3,3,4,4) satisfying Proposition 1 have to be the following two graphs.

For $G_1$, consider the pair of vertices $(v_1, v_7)$, the only 2-path between them is $v_1 v_5 v_7$, thus let $c(v_1 v_5) = 1, c(v_5 v_7) = 2,$. Similarly, $c(v_6 v_7) = 1, c(v_2 v_6) = c(v_3 v_6) = c(v_4 v_6) = 2$, $c(v_2 v_1) = c(v_3 v_1) = c(v_4 v_1) = 2$. Therefore, there is no rainbow $v_2 - v_3$ path.

For $G_2$, consider the pairs of vertices $(v_2, v_7), (v_3, v_7)$, we have $c(v_2 v_6) = c(v_3 v_6)$, consider the pairs of vertices $(v_2, v_5), (v_3, v_5)$, then $c(v_2 v_1) = c(v_3 v_1)$, so there is no rainbow $v_2 - v_3$ path.

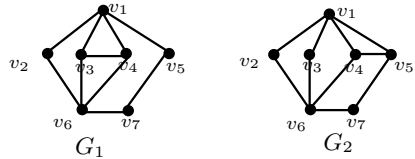

Figure 6: Graphs with degree sequence (2,2,2,3,3,4,4) satisfying Proposition 1.

(8) $\begin{cases} d_G = (4, 2, 2, 2, 2, 2, 2) \\ d_{\overline{G}} = (2, 4, 4, 4, 4, 4, 4). \end{cases}$

There is no graph with degree sequence (4,2,2,2,2,2,2) satisfying Proposition 1.

(9) $\begin{cases} d_G = (4, 4, 2, 2, 2, 2, 2) \\ d_{\overline{G}} = (2, 2, 4, 4, 4, 4, 4). \end{cases}$

The graph $G$ with degree sequence (4,4,2,2,2,2,2) satisfying Proposition 1 is the subgraph $G'$ of $G_1$ depicted in Figure 6 by deleting the edge $v_3v_4$. Consider the pairs of vertices $(v_2, v_7), (v_3, v_7)$, we have $c(v_2v_6) = c(v_3v_6)$, consider the pairs of vertices $(v_2, v_5), (v_3, v_5)$, we get $c(v_2v_1) = c(v_3v_1)$. Then there is no rainbow $v_2 - v_3$ path. Therefore $rc(G') \neq 2$.

$$(10) \begin{cases} d_G = (4, 4, 4, 2, 2, 2, 2) \\ d_{\overline{G}} = (2, 2, 2, 4, 4, 4, 4). \end{cases}$$

There is no graph with degree sequence (4,4,4,2,2,2,2) satisfying Proposition 1.

**Theorem 3** *For $n \geq 8$, the lower bound $rc(G) + rc(\overline{G}) \geq 4$ is best possible, that is, there are connected graphs $G$ and $\overline{G}$ on $n$ vertices, such that $rc(G) = rc(\overline{G}) = 2$.*

*Proof.* For $n = 8$, see figure $G_8$, colored with two colors, the heavy line with color 2, the others with color 1. It is easy to check that they are rainbow connected.

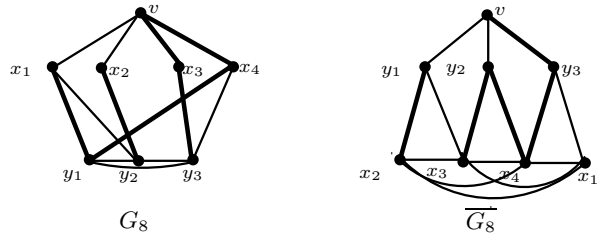

Figure 7: $rc(G) = rc(\overline{G}) = 2$ for $n = 8$.

If $n = 4k$, let $G$ be the graph with vertex set $X \cup Y \cup \{v\}$, where $X = (x_1, x_2, \cdots, x_{2k-1})$, $Y = \{y_1, y_2, \cdots, y_{2k}\}$, such that $N(v) = X$, $X$ is an independent set, $G[Y]$ is a clique, and for each $x_i$, $x_i$ is adjacent to $y_i, y_{i+1}, \cdots, y_{i+k}$, where the sum is taken modulo $2k$.

We define a coloring $c$ for the graph $G$ by the following rules:

$$c(e) = \begin{cases} 2 & \text{if } e = vx_i \text{ for } k+1 \leq i \leq 2k-1, \\ 2 & \text{if } e = x_iy_i \text{ for } 1 \leq i \leq 2k-1, \text{and } e = x_ky_{k+1}, \\ 1 & \text{otherwise.} \end{cases}$$

Then $c$ is a rainbow 2-coloring. And it is easy to check that $\overline{G}$ can also be colored by two colors and make it rainbow connected.

If $n = 4k + 1$, $G$ can be obtained by adding a vertex $x_{2k}$ to the vertex set $X$ in the case $n = 4k$, and joined $x_{2k}$ to $v, y_{2k}, y_1, \cdots, y_{k-1}$. With the coloring $c$ defined above, in addition with $c(vx_{2k}) = c(x_{2k}y_{2k}) = 2$, $G$ is rainbow connected.

If $n = 4k + 2$, $G$ can be obtained by adding two vertices $x_{2k}, y_{2k+1}$ to the vertex set $X$ and $Y$, respectively, in the case $n = 4k$, and joined $x_{2k}$ to $v$, $y_{2k+1}$ to each vertex in $Y$. And for each $x_i$, $x_i$ is adjacent to $y_i, y_{i+1}, \cdots, y_{i+k}$, where the sum is taken modulo $2k+1$. With the coloring $c$ defined above, in addition with $c(vx_{2k}) = c(x_{2k}y_{2k}) = 2$, $G$ is rainbow connected.

If $n = 4k+3$, $G$ can be obtained by adding two vertices $x_{2k+1}, y_{2k+1}$ to the vertex set $X$ and $Y$, respectively, in the case $n = 4k+1$, and joined $x_{2k+1}$ to $v$, $y_{2k+1}$ to each vertex in $Y$. And for each $x_i$, $x_i$ is adjacent to $y_i, y_{i+1}, \cdots, y_{i+k}$, where the sum is taken modulo $2k+1$, we also join $x_{k+1}$ to $y_{2k+1}$. With the coloring $c$ defined above, in addition with $c(vx_{2k+1}) = c(x_{2k+1}y_{2k+1}) = 2$, $G$ is rainbow connected.

**Theorem 4** *For $n = 4, 5$, $rc(G) + rc(\overline{G}) \geq 6$, and $rc(G) + rc(\overline{G}) \geq 5$ for $n = 6, 7$. All these bounds are best possible.*

*Proof.* From Theorem 2, we have $rc(G) + rc(\overline{G}) \geq 5$ for $n = 4, 5, 6, 7$.

For $n = 4$, as we have shown, $rc(G) + rc(\overline{G}) = 6$. If $n = 5$, the possible complementary connected graphs are:

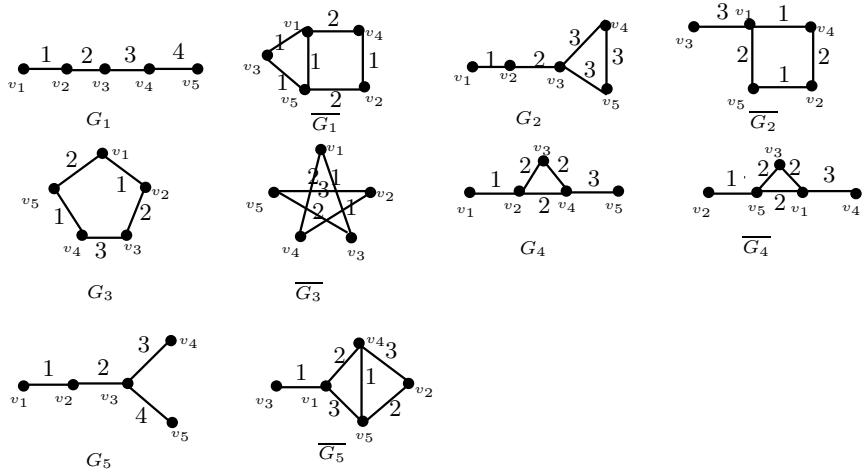

Figure 8: Complementary connected graphs for $n = 5$.

For all these cases, $rc(G) + rc(\overline{G}) \geq 6$.

For $n = 6$, let $G$ be the cycle $C_6$, whose vertices are $\{v_1, v_2, \cdots, v_6\}$. Then $rc(G) = 3$. We color the edges $v_1v_3$, $v_2v_4$, $v_3v_5$ in $\overline{G}$ by 2, and the other edges by 1. This coloring makes $\overline{G}$ rainbow connected. Therefore, $rc(G) + rc(\overline{G}) = 5$.

For $n = 7$, the graph $G_2$ in Figure 3 has rainbow connection number 2. We have shown that $rc(\overline{G_2}) \neq 2$, but we can use three colors to make it rainbow connected, just by assigning the edges $v_2v_4$ and $v_3v_5$ with color 3, the others the same as before. So, $rc(\overline{G_2}) = 3$. Thus, $rc(G) + rc(\overline{G}) = 5$.

## 4   Concluding remark

Given a graph $G$, a set $D \subseteq V(G)$ is called a *domination set* of $G$, if every vertex in $G$ is at a distance at most 1 from $D$. Further, if $D$ induces a connected subgraph of $G$, it is called a *connected dominating set* of $G$. The cardinality of a minimum connected dominating set in $G$ is called its *connected dominating number*, denoted by $\gamma_c(G)$. In [7], the authors proved that for every connected graph $G$ with minimum degree $\delta(G) \geq 2$, $rc(G) \leq \gamma_c(G) + 2$. In [5], the authors introduced a result of Nordhaus-Gaddum type result for the connected dominating number. They showed that if $G$ and $\overline{G}$ are both connected, then $\gamma_c(G) + \gamma_c(\overline{G}) \leq n+1$. If one uses their results, one can only get that $rc(G) + rc(\overline{G}) \leq \gamma_c(G) + \gamma_c(\overline{G}) + 4 \leq n+5$, which is weaker than our result.